%% file: main_file.tex
\RequirePackage[hyphens]{url}
\documentclass[review]{elsarticle}
\usepackage[margin=2.5cm]{geometry}
\usepackage{lineno,hyperref}

\usepackage{blindtext, graphicx}
\usepackage[cmex10]{amsmath}
\usepackage[cmex10]{mathtools}
\usepackage{amsfonts}
\usepackage{algorithmic}
\usepackage[font=small,labelfont=bf]{caption}
\usepackage[font=small,labelfont=bf]{subcaption}

\usepackage{bbm}
\usepackage{mathrsfs}

\DeclareMathOperator*{\argmin}{\arg\!\min}

\modulolinenumbers[5]

\makeatletter
\def\ps@pprintTitle{%
  \let\@oddhead\@empty
  \let\@evenhead\@empty
  \let\@oddfoot\@empty
  \let\@evenfoot\@oddfoot
}
\makeatother

\graphicspath{{./Pictures/}}

\begin{document}

\begin{frontmatter}

\title{Designing Near-Optimal Policies for Energy Management in a Stochastic Environment}

\input{author.tex}
\input{abstract.tex}
\input{els_keywords.tex}

\end{frontmatter}


\input{introduction.tex}
\input{model_description.tex}
\input{energy_management_problem.tex}
\input{energy_management_solutions.tex}
\input{results_discussion.tex}
\input{conclusion.tex}
\section*{Acknowledgments}
The authors would like to thank Cisco Systems, Inc. for its support.
\appendix
\input{appendix.tex}





\bibliography{references}

\end{document}

%% file: author.tex
\makeatletter
\def\@author#1{\g@addto@macro\elsauthors{\normalsize%
    \def\baselinestretch{1}%
    \upshape\authorsep#1\unskip\textsuperscript{%
      \ifx\@fnmark\@empty\else\unskip\sep\@fnmark\let\sep=,\fi
      \ifx\@corref\@empty\else\unskip\sep\@corref\let\sep=,\fi
      }%
    \def\authorsep{\unskip,\space}%
    \global\let\@fnmark\@empty
    \global\let\@corref\@empty  
    \global\let\sep\@empty}%
    \@eadauthor={#1}
}
\makeatother



\author{Chaitanya Poolla\corref{correspauthor}\fnref{chaitanya_affiliation}}
\ead{cpoolla@alumni.cmu.edu}
\fntext[chaitanya_affiliation]{Chaitanya Poolla is with Intel Corporation. This work was done while he was a graduate student at the department of Electrical and Computer Engineering, Carnegie Mellon University (SV).}

\author{Abraham K. Ishihara\fnref{abe_affiliation}}
\ead{abe.ishihara@west.cmu.edu}
\fntext[abe_affiliation]{Abraham K. Ishihara is with the department of Electrical and Computer Engineering, Carnegie Mellon University (SV).}

\author{Rodolfo Milito\fnref{rodolfo_affiliation}}
\ead{rmilito@starflownetworks.com}
\fntext[rodolfo_affiliation]{Rodolfo Milito is with Starflow Networks, Inc. This work was done while he was affiliated with Cisco Systems, Inc.}

\cortext[correspauthor]{Corresponding author}

%% file: abstract.tex
\begin{abstract}
With the rapid growth in renewable energy and battery storage technologies, there exists significant opportunity to improve energy efficiency and reduce costs through optimization. However, optimization algorithms must take into account the underlying dynamics and uncertainties of the various interconnected subsystems in order to fully realize this potential. To this end, we formulate and solve an energy management optimization problem as a Markov Decision Process (MDP) consisting of battery storage dynamics, a stochastic demand model, a stochastic solar generation model, and an electricity pricing scheme. The stochastic model for predicting solar generation is constructed based on weather forecast data from the National Oceanic and Atmospheric Administration. A near-optimal policy design is proposed via stochastic dynamic programming. Simulation results are presented in the context of storage and solar-integrated residential and commercial building environments. Results indicate that the near-optimal policy significantly reduces the operating costs compared to several heuristic alternatives. The proposed framework facilitates the design and evaluation of energy management policies with configurable demand-supply-storage parameters in the presence of weather-induced uncertainties.
\end{abstract}

%% file: els_keywords.tex
\begin{keyword}
Solar Energy; Energy Storage; Energy Management; Microgrid; Smart Building; Markov Decision Process
\end{keyword}

%% file: introduction.tex
\section{Introduction}
The electric grid was originally designed to support unidirectional power-flows from a few generating sources to a large number of consumers via transmission and distribution networks. With the recent growth in distributed renewable power generation the present grid is faced with the possibility of bidirectional power-flows. While greater renewable penetration is desirable from a sustainability perspective, the associated fluctuations pose significant grid stability challenges \cite{stadler2008power}. Further, growing concerns about efficiency, reliability, security, and carbon footprint necessitate the transformation of the present grid into a "smart grid" \cite{smartgrid}. In the proposed smart grid paradigm \cite{smartgrid}, information from various sensors are integrated via communication networks to enable intelligent real-time power-flow decisions. This results in several advantages including (i) reduced operation and maintenance costs via automation, (ii) increased efficiency by minimizing grid losses, (iii) increased renewable penetration via real time power management and integrated energy storage \cite{Rodrigues2014265}, and (iv) increased reliability. Thus, the smart grid is envisioned to integrate several forward-looking technologies while maintaining full backward compatibility with the existing grid without compromising grid stability \cite{farhangi2010}. The grid functions as a network of several interconnected nodes facilitating generation, transmission, and distribution of electric power at a regional scale or larger. A functional equivalent of the grid serving a smaller scale is known as a microgrid, which is regarded as a plug-and-play unit of the grid. Similarly, the functional equivalent of a microgrid within the scale of a single building is known as a nanogrid\footnote{In what follows, references to a microgrid or a nanogrid are made in the context of a smart grid indicating the presence of a communication network necessary for informed decision-making.} \cite{burmester2017review}.

The problem of making decisions for energy management is central to any electric grid including the microgrid or the nanogrid \cite{ZIA20181033}. Such decisions not only depend on grid-specific factors such as supply and demand levels but also depend on environment-specific factors such as the weather conditions, as described in the review \cite{AGUERAPEREZ2018265}. The study also highlights the indispensable nature of meteorological input for microgrid energy management systems, while noting the lack of importance provided to weather forecasts in several existing studies. In this study, we are concerned with decision problems in the presence of weather forecast-related uncertainties. In particular, we restrict our focus to a grid-connected nanogrid consisting of photovoltaic (PV) generation, battery, and load as shown in Figure \ref{FIG:nanogrid}. 

We discuss several previous works here and highlight key differences between them and the present work as applicable. A review of topologies and technologies for nanogrid energy management was presented in \cite{burmester2017review}. In \cite{zhang2014}, a distributed droop control for source converters combined with heuristic control laws for energy storage management is presented. On the other hand, we consider a non-controllable forecast-driven PV generation model with the objective of optimizing operational cost using energy storage. In \cite{zhang2016optimal}, a learning-based optimal demand response programming for home energy management was proposed. The study considered a static optimization problem subject to learning energy consumption models using neural networks and polynomial regression. By contrast, we consider a fixed demand model alongside battery dynamics in the optimization problem. In \cite{KRIETT2012residential}, a residential microgrid optimization problem with supply and demand management of electrical and thermal energy was considered. Several features including dishwasher, refrigerator, heating, electric aggregates, PV, and Combined Heat and Power (CHP) units were explicitly modeled. The optimization problem was formulated as a Mixed Integer Linear Program (MILP) embedded within an Model Predictive Control (MPC) scheme. A similar MPC-based approach with deterministic controllable load and PV are extended to a distributed setting in \cite{wang2014dynamic}. Unlike these MPC-based deterministic approaches, we consider an MDP formulation in the presence of stochastic non-controllable load and PV aggregate models, managed centrally by a Decision-Making Unit (DMU) as described in Section \ref{SEC:EMP}. \cite{zhang2018stochastic} proposed a stochastic MPC (SMPC) consisting of scenario generation and reduction phases in order to overcome the limitations of deterministic MPC schemes. In that study, the renewable and load forecasting errors were modeled as Gaussian white noise processes. By comparison, the forecasting errors in this work are modeled based on empirical distributions without predefined distributional assumptions. In \cite{riffonneau2011}, an optimal power flow problem for grid connected photovoltaic (PV) systems with batteries was considered. The study incorporated point estimates of weather forecasts and employed a deterministic dynamic programming-based approach, which was found to outperform a rule-based approach. However, in this study we incorporate weather forecast distributions and formulate the optimization problem as a Markov Decision Process (MDP) that is solved via stochastic dynamic programming. Similar to the studies in \cite{riffonneau2011} and \cite{HAWKES2007711}, we first propose heuristic policies followed by the near-optimal policy, whose performance is compared to its heuristic counterparts. In \cite{AMROLLAHI201766}, the impact of demand response implementation on component size optimization is studied. The problem was formulated as a MILP by considering the hourly mean PV generation, and the hourly mean load with a Gaussian multiplier. However, here we incorporate the uncertainties in the PV generation and the load by modeling them as cyclostationary stochastic processes with respective hourly empirical distributions. An MDP-based game between several residential customers and a utility service provider was proposed in \cite{misra2013}. The objectives were to enable strategic appliance load scheduling for the residential customers and strategic pricing schemes for the utility service provider. On the contrary, we employ an MDP formulation with the objective of minimizing operating costs in the presence of weather forecast-related uncertainties. Further, unlike the load scheduling problem addressed in the previous study, we consider the load to be non-controllable and modeled as stochastic.

An MILP-based optimization approach for solving a unit commitment problem was proposed in \cite{borghetti2008milp}. The approach involved the reduction of a mixed integer nonlinear structure into an MILP-compliant structure to enable the use of MILP solvers. In \cite{koot2005energy}, a vehicular energy management optimization problem was formulated with a nonlinear structure, which was reduced to a convex Quadratic Programming (QP) structure. In both of these works, reducing the nonlinear structure compromises the representational accuracy of the original problem. As noted in \cite{riffonneau2011}, most previous studies either use approaches that require specialized solvers or employ structural approximations that compromise the original nonlinear problem. Further, most studies capable of incorporating nonlinearity lack a probabilistic framework to accommodate uncertainty. To the best of our knowledge, no study has so far proposed a weather forecast-integrated MDP framework for designing a nanogrid Energy Management System (EMS) with a comprehensive treatment for handling state and decision constraints under uncertainty. The advantage of using an MDP formulation is that it accommodates both nonlinear dynamics and nonlinear objectives, unlike the more commonly used linear or quadratic programming formulations. Further, an MDP-based formulation naturally lends itself to incorporate uncertainties in the form of probability distributions unlike most MPC-based formulations.

Our main contributions include (1) a weather forecast-integrated MDP formulation of a nanogrid EMS with uncertain demand-supply capable of accommodating nonlinear objectives and nonlinear dynamics, (2) heuristic policy\footnote{In this work, heuristic policies are also referred to as \emph{naive} policies. The reader is referred to Section \ref{SSEC:naive_policies} for details.} designs that account for the present and lookahead demand-supply scenarios, and (3) a near-optimal energy management policy using Stochastic Dynamic Programming (SDP) along with a comprehensive treatment for handling state and decision constraints under uncertainty.


The remainder of this paper is organized as follows: Section \ref{SEC:Model_Desc} introduces the models of the various components within the nanogrid. The energy management problem is posed in Section \ref{SEC:EMP}. Section \ref{SEC:Policy} proposes energy management solutions in the form of policies. Heuristic policies are formulated in Section \ref{SSEC:naive_policies}. An optimal\footnote{Though we seek an optimal policy, the solution methodology involves the use of approximation resulting in a near-optimal policy, as presented in Section \ref{SEC:OPC}. A detailed description of the near-optimal policy design via approximate dynamic programming can be found in \cite{bertsekas2005dynamic}.} policy sought within the MDP framework is described in Section \ref{SEC:OPC}. Simulation results for residential and commercial buildings are presented and discussed in Section \ref{SEC:Discussion} followed by concluding remarks in Section \ref{SEC:Conclusion}. The treatment of state and decision constraints in the SDP formulation is discussed in \ref{app:constraint_handling}.


%% file: model_description.tex
\section{Model Description}
\label{SEC:Model_Desc}
\begin{figure}
\centering
\includegraphics[width=0.4\textwidth]{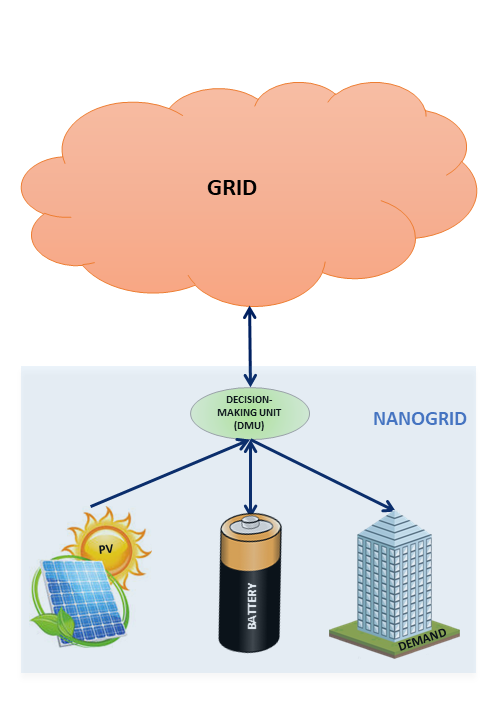}
\caption{Schematic of a grid-connected nanogrid}
\label{FIG:nanogrid}
\end{figure}

A grid-connected nanogrid or simply nanogrid, in this work, refers to a system consisting of a building load, PV generation, energy storage (battery), and network interconnections as depicted in Figure \ref{FIG:nanogrid}. In what follows, we describe the notation and the models of the individual components.
%
%
%
%
%
%
\subsection{Preliminaries}
\label{SSEC:Time}
Let the time interval of interest be denoted by $\mathcal{I} := [t_{st},t_{end}]$, where $t_{st}$ and $t_{end}$ denotes the start and end times, respectively. Let a partition of $\mathcal{I}$ be given by $\mathcal{T} := \{t_{k}\}_{k=0}^N$ such that the following properties hold true:
\begin{align}
\label{EQ:Time_Span}
t_0 &= t_{st} \nonumber\\
t_N &= t_{end}\nonumber \\
t_k &= t_0 + k\ \Delta t
\end{align}
where $\Delta t$ represents a uniform measure of subinterval $[t_{k-1},t_k]\ \forall\ k = \{1,2,\cdots,N\}$. 
For the purpose of integrating weather data from the past, let $\mathcal{I}_p$ denote a time interval containing elements which represent time in the past. Let a partition of $\mathcal{I}_p$ be given by $\mathcal{T}_p := \{T_{i}\}_{i=0}^{N_p}$ such that the following hold true:
\begin{equation}
\label{EQ:Time_Span_past}
T_i = T_0 + i\ \Delta T   
\end{equation}
where $\Delta T$ represents a uniform measure of subintervals $[T_{i-1},T_i]\ \forall\ i = \{1,2,\cdots,N_p\}$. For the decision problem described in Section \ref{SEC:Policy}, let each point $t_{k}$ in the sequence $\{t_{k}\}_{k=0}^N$ define a decision epoch, at which instant a decision-rule is formulated. This rule is designed to govern the decisions implemented during the interval $[t_k,t_{k+1})$. Further, let the measure of each subinterval be uniform such that $(T_i-T_{i-1})=\Delta T\ \forall i =\{1,\cdots,N_p\}$. In this work, we set $\Delta T = \Delta t = 3600$ seconds so that the time instants in $\mathcal{T}_p$ and $\mathcal{T}$  are one hour apart in time.
\subsection{Photovoltaic Generation}
\label{SSEC:PV}
Let the power output from the PV generation be represented by a discrete time stochastic process $e_k, \forall\ k \in \{0,\cdots,N\}$ with discrete states $pv_k \in \{e_k^{min},\cdots, e_k^{max}\}$. We assume this process to be cyclostationary with a time period of 24 hours. The probability distribution of the random variable $e_k$ at time instant $t_k \forall\ k \in \{0,\cdots,N\}$ is computed from the underlying irradiance distributions, which are inferred from datasets provided by the National Oceanic and Atmospheric Administration (NOAA). These irradiance distributions are converted into the corresponding PV output distributions, as described toward the end of Section \ref{NOAA}. We assume the PV output to be non-controllable in this work.
\subsubsection{Integrating Weather Forecasts from NOAA}
\label{NOAA}
NOAA provides a variety of weather forecast products varying in spatio-temporal resolution, prediction horizon, and update frequency \cite{noaa_ncep}. In this work, we use the forecast archives from NOAA's North American Mesoscale Forecast System (NAM) which provides forecasts extending up to 84 hours into the future. For sensor measurements, we use irradiance data from the pyranometers installed on Building 19 within the Carnegie Mellon University's (CMU's) campus located at Moffett Field, California. Alternatively, the NAM forecast data provides a 0-hr ahead forecast, which is regarded as an estimate of the true irradiance and may serve as a proxy if the sensor measurements are unavailable. Let $\mathcal{M}_i$ denote the measured irradiance value at a time instant in the past $T_i \in \mathcal{T}_p$ and let $\mathcal{F}_i$ denote the corresponding forecast value at $T_i$ based on a $h$-hour ahead forecast. Let the difference be represented by $\mathcal{E}_i$ as defined below:
\begin{equation}
\label{EQ:forecast_error}
\mathcal{E}_i = \mathcal{M}_i - \mathcal{F}_i
\end{equation}
The sensor data $\mathcal{M}_i$ is available at a sub-hourly resolution. However, the forecast data $\mathcal{F}_i$ from the NAM model is generated once every six hours at $\{00, 06, 12, 18\}$ hr UTC. Let the set of hours during which a forecast is generated be represented by $\mathcal{T}_f$. Owing to the differences in granularity between measurements and forecasts, we note that $\mathcal{E}_i$ from a $h$-hour ahead forecast is available only when $\big(\frac{T_i}{3600} \mod 23\big) = \big((F+h) \mod 23\big)\ \forall\ F \in \mathcal{T}_f$. However, for the purpose of inference, we assume that $\mathcal{E}_i$ is piecewise constant. Therefore, $\mathcal{E}_i$ is well-defined for every $i \in \{0,\cdots,N_p\}$. In this manner, the error values $\mathcal{E}_i$ are well-defined for every $T_i \in \mathcal{T}_p$.

Given the error population $\mathcal{E}_i\ \forall\ i \in \{0,\cdots,N_p\}$, let $\mathbf{E}_j$ denote the dataset available at the time instant $T_j \in \mathcal{T}_p$ containing the error data $\mathcal{E}_i$, where $i = \big\{m \in \{0,\cdots,N_p\} \mid (\frac{T_m}{3600} \mod 23) = (\frac{T_j}{3600} \mod 23)\big\}\ and\ T_m \leq T_j$. In other words, the dataset $\mathbf{E}_j$ is a collection of all the error data $\mathcal{E}_i$ such that the corresponding time instants $T_j$ and $T_i$ correspond to the same hour.
We now proceed to infer the error distributions. Let the error generating process be represented by a cyclostationary stochastic process $\hat{\mathcal{E}}_k\ \forall\ k \in \{0,\cdots,N\}$ with a time period of 24 hours. Let its distribution be denoted by $\mathscr{F}(D_k)$, where $D_k$ refers to the underlying dataset used for the inference. Under the periodicity assumption, we note that $\forall k\ \in \{0,\cdots,N\},\ \forall i\ \in \{0,\cdots,N_p\}$:
\begin{equation}
\label{EQ:past_dataset_inference}
D_k = E_i \Leftrightarrow \frac{t_k}{3600} = \frac{T_i}{3600} \\
\end{equation}

Given the characterization of the stochastic process $\hat{\mathcal{E}}_k, \forall\ k \in \{0, \cdots ,N\}$ as described above, we train a model to predict the corresponding measured value $\hat{\mathcal{M}}_k$ by leveraging the weather forecast at $t_k \in \mathcal{T}$. While any standard linear regression model such as the one in \cite{poolla2018solar} may be employed, we employ a signal-noise model as shown below.
\begin{equation}
\label{EQ:Signal_Noise_Model}
\hat{\mathcal{M}}_k = \mathcal{F}_k + \hat{\mathcal{E}}_k
\end{equation}
In this manner, the distribution of the estimated measurement $\hat{\mathcal{M}}_k$ is inferred with the knowledge of the forecast $\mathcal{F}_k$ and the error distribution $\hat{\mathcal{E}}_k$. Similarly, other variables representing different aspects of weather such as temperature, cloud cover, etc. can be estimated. With the knowledge of these weather fluctuation distributions, models that map the environment states to the solar power output can be used to obtain the solar power distribution. Although neural network models can be used \cite{poolla2014neural}, we use a linear model similar to the one described in \cite{sharma2016modeling}. Let the resulting solar power probability at $t_k \in \mathcal{T}$ be represented by $P^e(k,pv_k)$, where $pv_k \in \{e_k^{min}, \cdots, e_k^{max}\}$ denotes one of the states of the photovoltaic generation defined in Section \ref{SSEC:PV}.
\subsection{Load}
\label{SSEC:Load}
The demand in a residential or commercial building could be broadly categorized into Heating Ventilation and Air Conditioning (HVAC), lighting, appliance and plug loads. In this work, we employ a data-based load model for the aggregate demand in the nanogrid. The hourly load data was obtained from the openEI database \cite{openei_data} for both the residential and the commercial buildings. Similar to the PV generation model, the nanogrid load is modeled as a discrete time discrete state stochastic process $l_k, \forall\ k \in \{0, \cdots ,N\}$ with states $lo_k \in \{l^{min}, \cdots, l^{max}\}$. This process is assumed to be cyclostationary with a time period of 24 hours. Let the probability distributions of $l_k$ be denoted by $P^l(k,lo_k)$, where $lo_k \in \{l^{min}, \cdots, l^{max}\}$ denotes a state of the load process. In this work, we assume the demand cannot be controlled and is always fulfilled.
\subsection{Energy Storage}
\label{SSEC:Storage}
The nanogrid energy storage element (battery) is modeled as a dynamical system with capacity $\mathcal{S}$. Given the non-controllable nature of both the PV generation and the load, the battery offers to control the power flow in the grid-connected nanogrid (Equation \ref{EQ:power_balance}). While several dynamical models with varying levels of complexity have been considered \cite{he2011evaluation}, we use the following linear dynamical model as shown below.
\begin{equation}
\label{EQ:batt_dyn}
s_{k+1} = \eta_s s_k - \xi_p v_k \Delta t
\end{equation}
where, $s_k \in [0,\mathcal{S}]$ represents the state of the battery and $v_k$ represents the net power output from the battery, both at $t_k \in \mathcal{T}$.
The power flow constraints associated with the battery are represented by $v_k \in [P_{min},P_{max}]$. The handling of the state and power constraints are described in \ref{app:constraint_handling}. Further, let the losses due to capacity fading and power fading be represented by the corresponding efficiencies $\eta_s$ and $\xi_p$\footnote{We note that $\xi_p < 1.0$ during the charge phase and $\xi_p > 1.0$ during the discharge phase. In other words, more power is discharged from the battery than what is output and less power is used to charge the battery than what is input.}, respectively \cite{vetter2005ageing}. 
\subsection{Grid Transactions}
\label{SSEC:grid}
The nanogrid transacts power with the grid at an associated cost. In this work, we consider these transactions to be lossless and unconstrained. The decisions driving these transactions are enabled by the Decision-Making Unit (DMU). A functional representation of the nanogrid along with its DMU is shown in Figure \ref{FIG:nanogrid}. 
Let the power flowing from the grid to the DMU at any $t_k \in \mathcal{T}$ be denoted by $u_k$\footnote{By convention, we regard the power flowing into the DMU as positive.}. The transaction cost per unit of power is determined by the utility pricing scheme. We consider deterministic pricing schemes for the residential and commercial building scenarios. For the purchase pricing scheme ($c_p(k), k \in \{0, \cdots, N\}$), we use the $PG\&E$ schedule E6 \cite{pge_residential_e6} for the residential building and the $PG\&E$ schedule A10 \cite{pge_commercial_a10} for the commercial building. For the selling pricing scheme ($c_s(k), k \in \{0, \cdots, N\}$), we use the $PG\&E$ E-SRG \cite{pge_srg_ppa} schedule for both the residential and commercial buildings. Thus, the monetary cost incurred due to the grid transactions over the horizon $\mathcal{T}$ can be written as:
\begin{equation}
\label{EQ:cost_fcn}
\mathbf{C}(u_0,\cdots,u_{N-1};w_0,\cdots,w_{N-1}) = \sum_{i=0}^{N-1} \{c_p(i)I(u_i>0) + c_s(i)I(u_i<=0)\}u_i \Delta t
\end{equation}
where, a positive value for $\mathbf{C}$ indicates the cost to be paid by the nanogrid to the grid on the account of the grid transactions $(u_0,\cdots,u_{N-1})$.
\subsection{System Parameters}
\label{SSEC:sys_params}
Let the set of system parameters at $t_k \in \mathcal{T}$ be denoted by $\mathcal{W}_k$. Of these, let the set of deterministic parameters be denoted by $z_k$ and let the set of distributions of the corresponding random variables be denoted by $w_k$. The deterministic parameters consist of $z_k = [c_p(t_k),c_s(t_k),\mathcal{S},P_{min},P_{max},\eta_s,\xi_p,N,\Delta t]$ and the distributions of the random variables $(e_k,l_k)$ consist of $w_k = [P^e(k,\cdot),P^l(k,\cdot)]$.

%% file: energy_management_problem.tex
\section{Energy Management Problem}
\label{SEC:EMP}
The energy management problem in the nanogrid amounts to determining the power flows $u_k$ and $v_k$ at each time $t_k \in \mathcal{T}$. Candidate solutions to this problem must satisfy the state and power constraints described in Equations \ref{EQ:batt_dyn}, \ref{EQ:PC2CC}, and \ref{EQ:power_balance}.  Let the PV output and the load at time $t_k \in \mathcal{T}$ be denoted by $e_k^{(r_e)}$ and $l_k^{(r_l)}$, respectively. Here $e_k^{(r_e)}$ and $l_k^{(r_l)}$ represent the realizations of the underlying stochastic processes $e_k$ and $l_k$, respectively. Given the values of the PV output and the load, the decisions computed at the DMU ($u_k$, $v_k$) must result in a power balance at each time $t_k \in \mathcal{T}$, as expressed in \ref{EQ:power_balance}.
\begin{equation}
\label{EQ:power_balance}
e_k^{(r_e)} + u_k + v_k + l_k^{(r_l)} = 0\ 
\end{equation}
\subsection{Optimal Energy Management Problem}
\label{SEC:OPF}
We are interested in the minimum cost solution to the energy management problem. We consider a metric that measures the expected monetary cost over the finite horizon $\mathcal{I}$. This metric is a function of both the grid transactions $u_k,\ \forall k \in \{0, \cdots, N-1\}$ and the expected storage state at the end of the horizon. We note that the storage state at the end of the horizon $t_N$ is unrealized at every $t_k < t_N$, and is regarded as a random variable. Let this random variable be denoted by $\hat{s}_N$. The probability distribution of $\hat{s}_N$ can be computed with the knowledge of (i) the PV and the load distributions $P^e(i,pv_i)$ and $P^l(i,lo_i)$, respectively, over all instants $i \in \{k,\cdots,N-1\}$, and (ii) the storage dynamics described in Equation \ref{EQ:batt_dyn}. Therefore, the cost metric over the horizon $[t_k,t_N]$, also known as the cost-to-go from stage $k$ to stage $N$, can be written as:
\begin{multline}
\label{EQ:Cost}
J_1(s_k,u_k,u_{k+1},\cdots,u_{N-1};\mathcal{W}_k,\cdots,\mathcal{W}_N) = \mathbb{E} \bigg\{\sum_{i=k}^N L(u_i;\mathcal{W}_i) + g_1(\hat{s}_N;\mathcal{W}_k,\cdots\mathcal{W}_N)\bigg\} \\
= \mathbb{E} \bigg\{\sum_{i=k}^N \big(\mathbbm{1}(u_i>0)c_p(t_i) + \mathbbm{1}(u_i\leq0)c_s(t_i)\big)u_i - \hat{s}_N c_s(t_N)\bigg\}
\end{multline}
where, $\mathbbm{1}(\cdot)$ denotes the indicator function, $\mathbb{E}$ denotes the expectation operator with respect to the PV and load distributions, and $\mathcal{W}_k$ denotes the system parameters at time $t_k \in \mathcal{T}$. While the above metric captures the essential monetary costs, we define another cost metric to emphasize upon a full battery reserve at the end of the horizon. This is accomplished by introducing a multiplier $m>>1$, which adjusts the terminal cost to $\mathbb{E}\{g(\hat{s}_N;\mathcal{W}_k,\cdots,\mathcal{W}_N)\} = m(\mathcal{S}-\mathbb{E}\{\hat{s}_N\})c_s(t_N)$. Thus, the new cost metric over the horizon $[t_k,t_N]$ can be written as:
\begin{equation}
\label{EQ:cost_new}
J(s_k,u_k,u_{k+1},\cdots,u_{N-1};\mathcal{W}_k,\cdots,\mathcal{W}_N) = \mathbb{E} \bigg\{\sum_{i=k}^N L(u_i;\mathcal{W}_i) + g(\hat{s}_N;\mathcal{W}_k,\cdots\mathcal{W}_N)\bigg\}
\end{equation}
We describe the solution to the optimal energy management problem in Section \ref{SEC:OPC}.

%% file: energy_management_solutions.tex
\section{Energy Management Solutions}
\label{SEC:Policy}
The solution to the energy management problem described in Section \ref{SEC:EMP} are presented here. The power flow decisions ($u_k,v_k$) can be viewed as the result of a mapping from the state space to the decision space of the system. Such a mapping is also known as a decision-rule or policy. Let the nanogrid energy management policy be denoted by $\pi$ and be comprised of a grid transaction policy $\pi_G$ and a storage transaction policy $\pi_B$ which determine the power flows $u_k$ and $v_k$, respectively. Thus, the policy $\pi$ can be written as:
\begin{eqnarray}
\label{EQ:policy}
\pi \coloneqq (\pi_G,\pi_B), \mbox{where} \notag \\
\pi_B: \{0, \cdots, N-1\} \times \mathcal{S}_k \to \mathcal{V}_{s_k} \notag \\
\pi_G: \{0, \cdots, N-1\} \times \mathcal{S}_k \to \mathcal{U}_{s_k}
\end{eqnarray}
where, the admissible decision space is represented by $(\mathcal{U}_{s_k},\mathcal{V}_{s_k}) \subset \mathbf{R}^2$ and $\mathcal{S}_k$ denotes the storage state space. It must be noted that the admissible decision space consists of the decisions which, when implemented do not result in the violation of the system constraints (refer to \ref{app:constraint_handling} for details). Also, from Equation \ref{EQ:power_balance}, we note that $u_k$ and $v_k$ are not independently determined. Therefore, the corresponding policies $\pi_G$ and $\pi_B$ are dependent. In addition, we note that the performance of various policies can be compared based on cost metrics such as the one in Equation \ref{EQ:cost_fcn}. For purposes of benchmarking, we first describe heuristics-based naive policies. This is followed by the description of the near-optimal policy.
\subsection{Naive Policies}
\label{SSEC:naive_policies}
\begin{figure*}
\centering
\includegraphics[width=\textwidth]{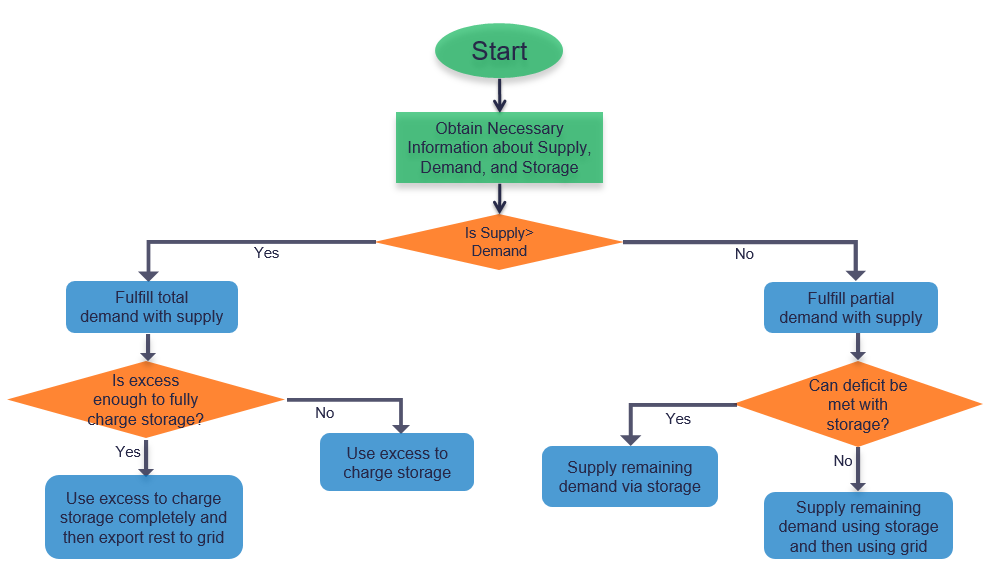}
\caption{Flow Chart depicting the Exhaustive Storage Dependence Policy}
\label{FIG:policy_1}
\end{figure*}

\begin{figure*}
\centering
\includegraphics[width=\textwidth]{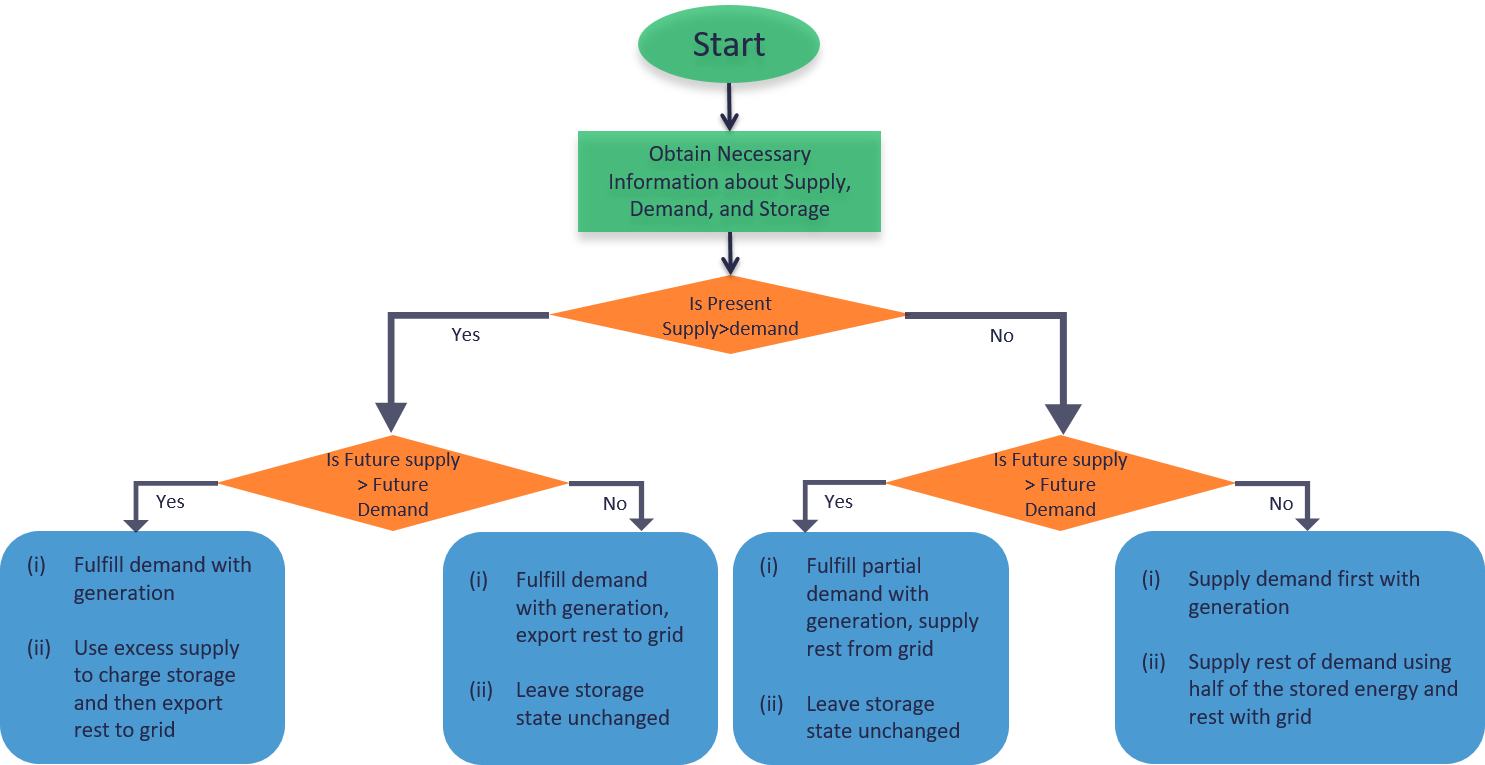}
\caption{Flow Chart depicting the three-hour Lookahead Policy}
\label{FIG:policy_2}
\end{figure*}

We employ heuristics based on demand-supply characteristics in the present and in the expected future to formulate the naive policies. Each naive policy along with its decision-making considerations is presented below.
\begin{enumerate}
\item \emph{Policy 1: Exhaustive Storage Dependence Policy} ($\pi_1$): Given the non-controllable supply $e_k^{(r_e)}$ and demand $l_k^{(r_l)}$, this policy is designed to bridge the demand-supply gap by depending on the storage transactions $v_k$. When the storage resource is no longer usable, the policy resorts to the grid transactions $u_k$ to bridge the demand-supply gap. The decision-making process for this policy ($\pi_1$) is illustrated by the flowchart in Figure \ref{FIG:policy_1}.
\item \emph{Policy 2: Finite Horizon Lookahead Policy with a three-hour lookahead} ($\pi_2$): This policy is designed to make informed decisions based on the present nanogrid state as well as the expected state over a finite horizon (future). We consider a three hour ahead horizon in this work. With these considerations, four scenarios are possible: (i) excess supply in the present and excess supply cumulatively expected over the horizon, (ii) excess supply in the present and deficit supply cumulatively expected over the horizon, (iii) deficit supply in the present and excess supply cumulatively expected over the horizon, and (iv) deficit supply in the present and deficit supply cumulatively expected over the horizon. When there is excess supply in the present and excess supply is cumulatively expected over the horizon, the storage resource is charged before depending on the grid for power balance, thereby accommodating excess generation during the present and over the future horizon relying primarily on the storage resource. In the scenario with deficit supply in the present and deficit supply cumulatively expected over the horizon, the storage resource is half-discharged to meet the deficit before depending on the grid. The rationale behind this mechanism is that the storage is half-used to meet the present deficit and the rest is retained to meet the deficit expected over the future horizon. In the remaining scenarios, the policy relies on the generation and the grid for power balance, leaving the storage resource unused. The decision flow corresponding to the policy $\pi_2$ is illustrated by the flowchart in Figure \ref{FIG:policy_2}.
\end{enumerate}
\input{optimal_policy_computation.tex}

%% file: optimal_policy_computation.tex
\subsection{Near-Optimal Policy Computation}
\label{SEC:OPC}
We address the optimal energy management problem posed in Section \ref{SEC:OPF} here. First, we note that the states ($s_k$), the decisions $(u_k,v_k)$, the Markovian state transitions (Equation \ref{EQ:batt_dyn}), and the cost metric $J$\footnote{Though the optimal policy minimizes the metric $J$, we evaluate the policies in Section \ref{SEC:Discussion} based on their actual monetary cost $J_1$.} (Equation \ref{EQ:cost_new}) together constitute a Markov decision process (MDP). Within this framework, we note that a decision $(u_k,v_k)$ enforced on the system at $t_k \in \mathcal{T}$ results in the state randomly transitioning from $s_k$ at $t_k$ to $s_{k+1}$ at $t_{k+1}$. The cost incurred during this transaction is represented by $L(u_k;\mathcal{W}_k)$ as shown in Equation \ref{EQ:cost_new}. 
%

Let the optimal policy be denoted by $\pi^*$. The objective of $\pi^*$ is to minimize the expected cost over the optimization horizon $[t_0,t_N]$. Accordingly, we pose the optimization problem below:
\begin{eqnarray}
\label{EQ:opt_policy}
\pi^*(s_0;(\mathcal{W}_k)_{k=0}^N) = (u^*_k,v^*_k)_{k=0}^{N-1},\ \mbox{such that} \notag \\
(u^*_k)_{k=0}^{N-1} = \argmin_{\big[(u_k)_{k=0}^{N-1} \in \prod\limits_{k=1}^{N-1} \mathcal{U}_{\mathbb{E}\big\{\hat{s}_k \mid s_0\big\}}\big]} J(s_0,u_0,\cdots,u_{N-1},\{\mathcal{W}_k\}_{k=0}^N), \mbox{and} \notag \\
v^*_k = -(u^*_k+\mathbb{E}(e_k)+\mathbb{E}(l_k)),\ \forall k \in \{0,\cdots,N-1\}
\end{eqnarray}
where, $\hat{s}_k$ denotes a random variable that represents the uncertain storage state at a future instant $t_k \in \mathcal{T}$. The corresponding optimal cost then becomes:
\begin{equation}
\label{EQ:opt_cost}
V_0(s_0;\{\mathcal{W}_k\}_{k=0}^N) = J(s_0,u^*_0,\cdots,u^*_{N-1};\{\mathcal{W}_k\}_{k=0}^N)
)\end{equation}
where, $V_0(s_0;\{\mathcal{W}_k\}_{k=0}^N)$ denotes the value of the state $s_0$ at the time instant $t_0$. The value function $V$ maps the states and parameters $(s_k;\{\mathcal{W}_k\}_{k=0}^N)$ at $t_k$ to their minimum cost as shown below:
\begin{eqnarray}
\label{EQ:VF1}
V_k(s_k;\{\mathcal{W}_i\}_{i=k}^N) = J(s_k,u^*_k,\cdots,u^*_{N-1};\{\mathcal{W}_i\}_{i=k}^N)
\end{eqnarray}
In other words, the value function at the instant $k$ describes the optimal cost-to-go from stage $k$ through the final stage $N$.

We address the above discrete time stochastic dynamic optimization problem by applying the principle of optimality and solving the resulting sub-problems using a dynamic programming (DP) formulation. The solutions to the DP formulation are then obtained by backward induction. In what follows, we describe the DP formulation of the problem presented in Equation \ref{EQ:opt_policy} and seek the corresponding optimal solution.
By applying the principle of optimality, we can transform the optimal decision problem into a sequence of sub-problems as shown:
\begin{equation}
\label{EQ:principle_of_optimality}
J(s_0,u^*_0\cdots,u^*_{N-1};\mathcal{W}_0,\cdots,\mathcal{W}_N) = \min_{ u_0 } \Bigg\{ L(u_0;\mathcal{W}_0)\ + \left[ \sum_{k=1}^{N-1} \mathbb{E}\Big\{J\big(\hat{s}_k,u^*_k,\cdots,u^*_{N-1};\mathcal{W}_{k-1},\cdots,\mathcal{W}_N\big) \Big\}\right] \Bigg\}
\end{equation}
Using the definition of the value function from Equation \ref{EQ:VF1}, we rewrite Equation \ref{EQ:principle_of_optimality} as:
\begin{equation}
\label{EQ:poo_eqn}
V_0(s_0;\{\mathcal{W}_k\}_{k=0}^N) = \min_{ u_0 } \Bigg\{ L(u_0;\mathcal{W}_0)\ + \mathbb{E}\Big\{V_1(\hat{s}_1;\{\mathcal{W}_k\}_{k=0}^N) \Big\} \Bigg\}
\end{equation}
Continuing the course of breaking down into sub-problems, we arrive at the following generalized equation:
\begin{equation}
\label{EQ:Bellman_eqn}
V_k(s_k;\{\mathcal{W}_i\}_{i=k}^N) = \min_{ u_k } \Bigg\{ L(u_k;\mathcal{W}_k)\ + \mathbb{E}\Big\{V_{k+1}(\hat{s}_{k+1};\{\mathcal{W}_i\}_{i=k}^N) \Big\} \Bigg\}
\end{equation}
which is known as the Bellman equation, a recursive equation to update values at $t_k$ based on the values at the next time step $t_{k+1}$. We apply backward induction to solve for the values in the Bellman equation using a numerical approach. First, we observe that the domain of the value function is the continuous state space $s_k \in [0,\mathcal{S}]$. However, for purposes of numerically computing the values and the optimal decisions via the Bellman equation, we quantize the state space \footnote{By quantizing the state space for value computation, we approximate the value function defined over the continuous state space $[0,\mathcal{S}]$.}, thereby introducing approximation into the solution. 

Let the quantized state space be represented by the finite sequence $\mathbf{S} \coloneqq\ \{\mathbf{s}_1=0,\mathbf{s}_2,\cdots,\mathbf{s}_{N_s}=\mathcal{S}\}$. At the end of the horizon $t_N$, the values are computed for the quantized state space $s_N \in \mathbf{S}$:
\begin{equation}
V_N(s_N;\mathcal{W}_N) = g(s_N;\mathcal{W}_N)
\end{equation}
At every other time instant in the sequence $\{t_k\}_{k=N-1}^0$, the decisions resulting in the minimum cost-to-go are computed by solving the Bellman equation (Equation \ref{EQ:Bellman_eqn}) across the feasible decision space $\mathcal{U}_{s_k}$ (see \ref{app:constraint_handling}). It must be noted that we obtain a sequence of near-optimal decisions $\{u^*_k\}_{k=0}^{N-1}$ on account of quantizing the state space and approximating the value function. Thus, we obtain a near-optimal policy which approximates the optimal policy $\pi^*$.

%% file: results_discussion.tex
\section{Results and Discussion}
\label{SEC:Discussion}
The results from the NAM forecasts, the cyclostationary PV and load models, and the nanogrid simulations under the action of the various policies are described below.

\subsection{Weather Results}
\label{weather_results}
\begin{figure}
\centering
\begin{minipage}[t]{0.5\textwidth}
	\centering\captionsetup{width=.8\linewidth}%
	\includegraphics[width=\linewidth]{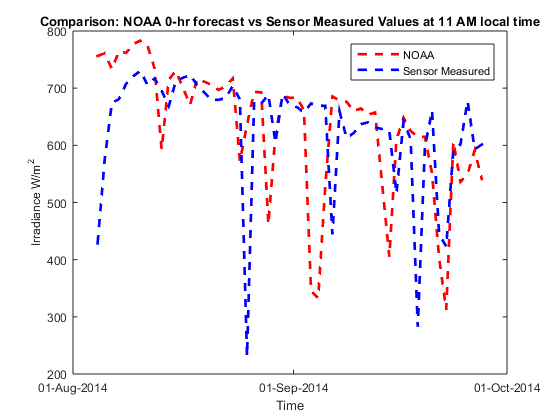}
	\caption{NOAA and Sensor irradiance data comparison at 11 AM, RMS Accuracy = $82.2\%$}
	\label{fig:11am_noaa_sensor}
\end{minipage}\hfill
\begin{minipage}[t]{0.5\textwidth}
	\centering\captionsetup{width=.8\linewidth}%
	\includegraphics[width=\linewidth]{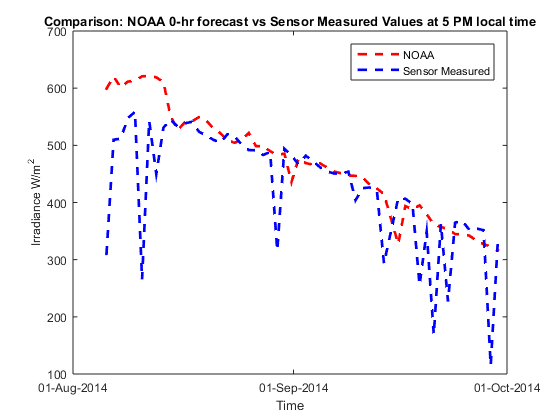}
	\caption{NOAA and Sensor irradiance data comparison at 5 PM, RMS Accuracy = $81.8\%$}
	\label{fig:5pm_noaa_sensor}
\end{minipage}
\end{figure}

\begin{figure}
\centering
\begin{minipage}[t]{0.5\textwidth}
	\centering\captionsetup{width=0.9\linewidth}%
	\includegraphics[width=\linewidth]{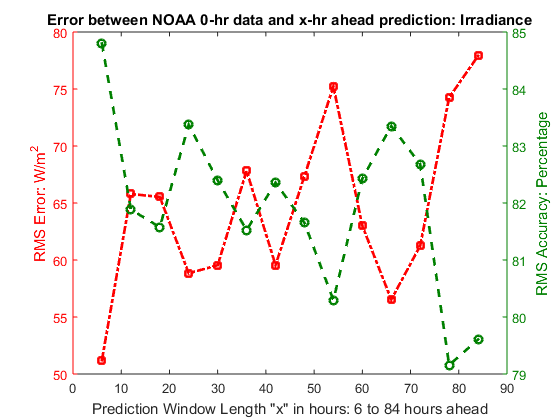}
    \caption{Comparison between NAM truth and $x$-hour ahead forecasts: $x \in \{6,12,\cdots,84\}$}
    \label{fig:hour_ahead_errors}
\end{minipage}\hfill
\begin{minipage}[t]{0.5\textwidth}
	\centering\captionsetup{width=0.9\linewidth}%
	\includegraphics[width=\linewidth]{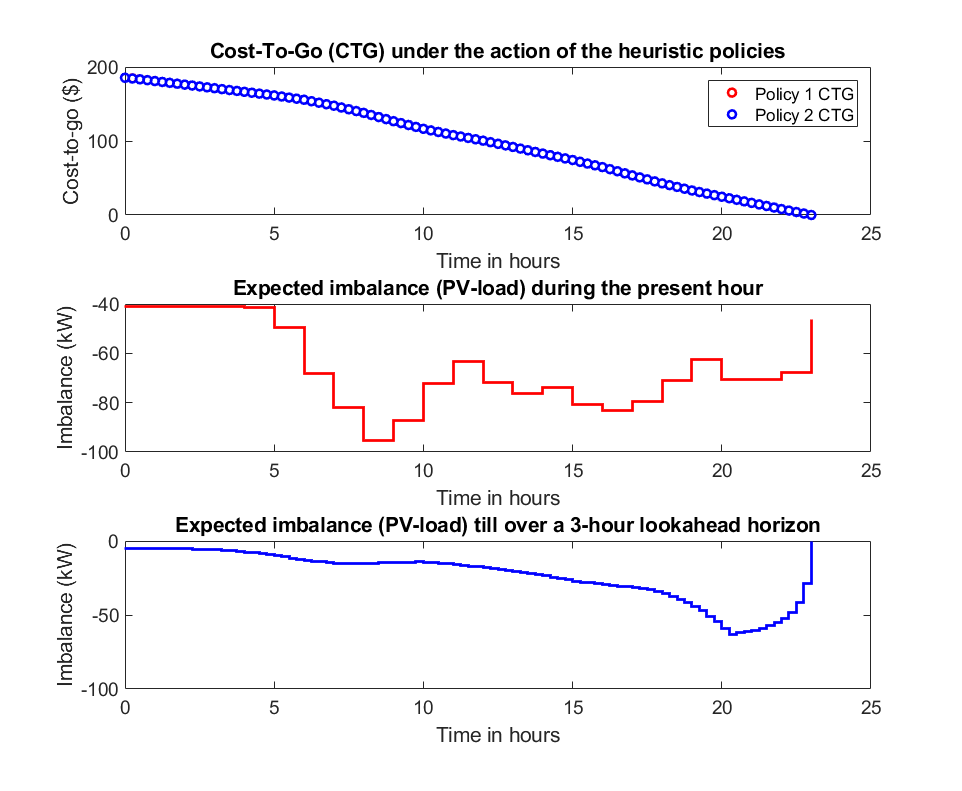}
	\caption{Cost-to-go under the action of the naive policies and the demand-supply imbalance scenarios in the commercial building simulation}
	\label{FIG:commercial_heuristic_ctg}
\end{minipage}
\end{figure}

We analyzed the NAM forecast data in relation to the sensor measurements and the NAM 0-hr ahead data (NAM truth estimates). First, we examined the similarity between the NAM truth estimates and the sensor measurements over a period of two months (Aug - Oct 2014). Comparisons at 11 PM PT and 5 PM indicate that the sensor data and the NAM truth estimates exhibit similar trends as shown in Figures \ref{fig:11am_noaa_sensor} and \ref{fig:5pm_noaa_sensor}. However, there exist considerable difference among both datasets. The Root Mean Square (RMS) error percentages computed indicated a mean error of $17.9\%$ along with a standard deviation of $33.4\%$. These error characteristics can be attributed to several factors including coarse spatial granularity, lack of local information about shadows or dust patterns, modeling error, and low update frequency of the NAM model (four per day).
Second, we compared the NAM forecasts to the NAM truth estimates. The NAM model provides forecasts up to 84 hours ahead. We present both the error and accuracy behavior across varying forecast horizon lengths in Figure \ref{fig:hour_ahead_errors}. The left vertical axis (in red) represents the RMS error scale and the right vertical axis (in green) represents the percentage RMS accuracy scale. The percentage accuracies are computed relative to the NAM truth estimates. It was found that the accuracy across various prediction horizons was consistent with a mean of $81.9\%$ and a standard deviation of $1.5\%$. In other words, the difference in forecast horizon length did not result in significant deviations from the mean percentage accuracy.
Third, we obtained the forecast error distributions from the weather archives as mentioned in Section \ref{NOAA}. These distributions were used to infer the distributions of the measured irradiance with the knowledge of the forecasts based on Equation \ref{EQ:Signal_Noise_Model}.
\subsection{Supply and Demand Model results}
\begin{figure}
\centering
\begin{minipage}[t]{0.5\textwidth}
	\centering\captionsetup{width=.9\linewidth}%
	\includegraphics[width=\linewidth]{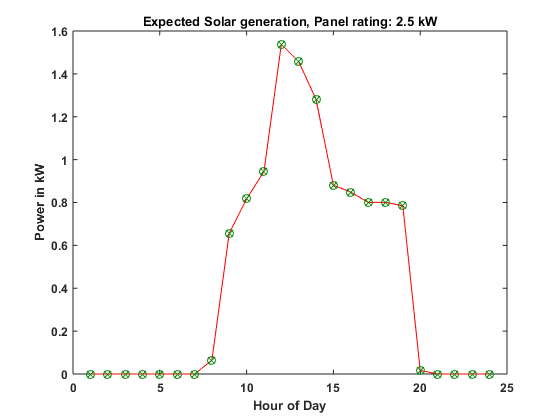}
	\caption{Expected Residential Solar Generation based on NOAA NAM irradiance estimates, $N_{pv}=5$}
	\label{FIG:solar}
\end{minipage}\hfill
\begin{minipage}[t]{0.5\textwidth}
	\centering\captionsetup{width=.9\linewidth}%
	\includegraphics[width=\linewidth]{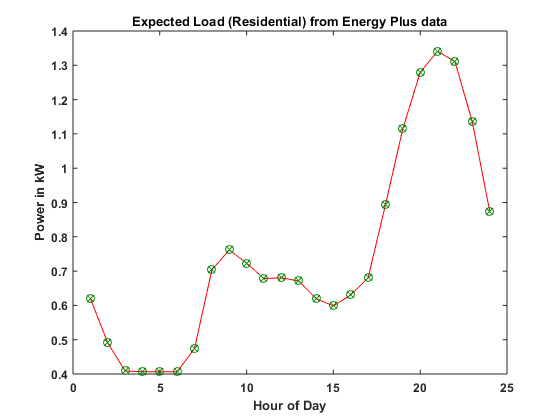}
	\caption{Expected Hourly Residential Load at San Francisco \cite{openei_data}, $N_{lo}=5$}
	\label{FIG:load}
\end{minipage}
\end{figure}

Based on the measured irradiance distributions and a linear model similar to the one described in \cite{sharma2016modeling}, the distributions of the cyclostationary PV generation process were inferred. For the residential PV installation, a capacity of $2.5 kW$ under standard conditions was used, and the expected PV generation is shown in Figure \ref{FIG:solar}. It can be observed that non-zero PV generation is expected over half of the 24 hours. In case of the commercial PV installation, a capacity of $84 kW$ under standard conditions was used based on the system capacity recommendations from \cite{davidson2015nationwide}.
Similar to the PV distributions, the distributions of the cyclostationary load process were inferred from yearly building power consumption data provided by OpenEI \cite{openei_data}. For the residential load scenario, the dataset corresponding to a residential building in San Francisco was used. For the commercial load scenario, the dataset corresponding to a commercial medium office building in Moffett Field was used. The hourly expected load corresponding to the residential building is depicted in Figure \ref{FIG:load}. While the load variation appears relatively smooth due to coarse temporal granularity, it must be noted that the daily load curves representing building demand vary abruptly in real time \cite{mathieu2011quantifying}.

\subsection{Nanogrid Simulation}
We simulated the nanogrid model under the action of both the naive and the near-optimal policies. The parameters of the battery were obtained from the specification sheet in  \cite{motors2015powerwall}. For the sake of analysis, the battery efficiencies were set to $\eta_s = \xi_p = 1.0$. The results are presented for two optimization horizons. First, a daily (24-hour) horizon was chosen to analyze the action of the naive and the near-optimal policies on the evolution of the battery state. This analysis was carried out for the residential building scenario. Second, a monthly (30-day) horizon was chosen to evaluate the monthly operating costs incurred by each policy. This evaluation was performed for both the residential and the commercial building scenarios.

\subsubsection{Residential building simulation over a daily horizon}
\begin{figure}
\centering
\begin{minipage}[t]{0.5\textwidth}
	\centering\captionsetup{width=.9\linewidth}%
	\includegraphics[width=\linewidth]{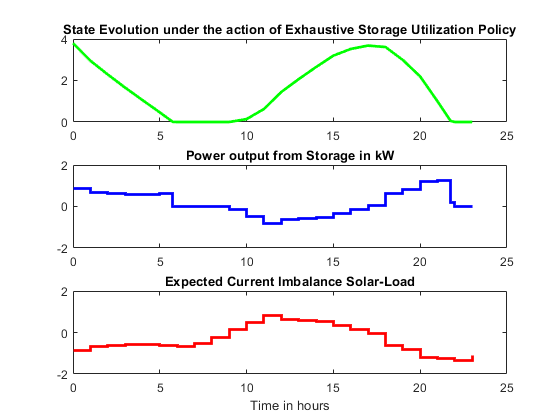}
	\caption{System under the action of Policy 1}
	\label{fig:policy1}
\end{minipage}\hfill
\begin{minipage}[t]{0.5\textwidth}
	\centering\captionsetup{width=.9\linewidth}%
	\includegraphics[width=\linewidth]{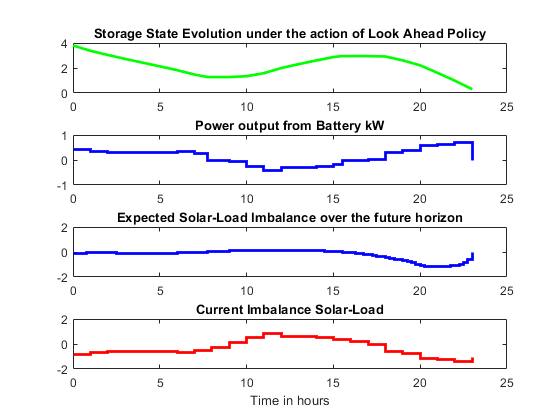}
	\caption{System under the action of Policy 2}
	\label{fig:policy2}
\end{minipage}
\end{figure}

The expected state evolution under the action of policy $\pi_1$ is shown in Figure \ref{fig:policy1}. The policy $\pi_1$ was designed to primarily depend on the battery for power balance prior to depending on the grid. In other words, we expect $\pi_1$ to charge the battery during excess PV generation and discharge the battery during deficit PV generation. When the battery cannot be further charged or discharged, we expect the policy to utilize the grid for power balance. The same behavior can be observed from the bottom subplot of Figure \ref{fig:policy1}. Deficit PV generation was present during the $[0,9]\cup[17,23]$ hour interval and therefore the policy discharged the battery. This is visible from subplot showing positive power output from the battery during these hours. In the case of excess generation during the $[10,16]$ hour interval, $\pi_1$ charges the battery as expected.
In case of the finite horizon look ahead policy $\pi_2$, the decisions are influenced by both the current demand-supply offset and the expected offset over a three-hour horizon into the future. The corresponding decision-making process is shown in Figure \ref{FIG:policy_2}. The action of the policy on the state evolution and the decision are shown in Figure \ref{fig:policy2}. To analyze the results, we observe three scenarios in relation to the sign of the imbalance (PV-load) between the present and expected future. The first scenario consists of the $[9,15]$ hour interval during which excess generation was observed in the present and expected in the future. In this case, as expected, the look ahead policy charged the battery as evident from the \emph{Power output from Battery kW} subplot in Figure \ref{fig:policy2}. The second scenario arises during the $[0,8)\cup(17,23]$ hour interval during which a deficit is observed in the present and expected in the future. As expected in this case, the policy discharges the battery. This behavior is evident from the subplot titled, \emph{Power output from Battery kW}. The third scenario arises during the hour interval $[8,9]\cup(15,17]$ when the sign of the present imbalance differs from the sign of expected future imbalance. In this case, the policy is designed to not alter the battery state. This behavior can be observed from the same subplot in Figure \ref{fig:policy2}.

\begin{figure}
\centering
\begin{minipage}[t]{0.5\textwidth}
	\centering\captionsetup{width=0.9\linewidth}%
	\includegraphics[width=\linewidth]{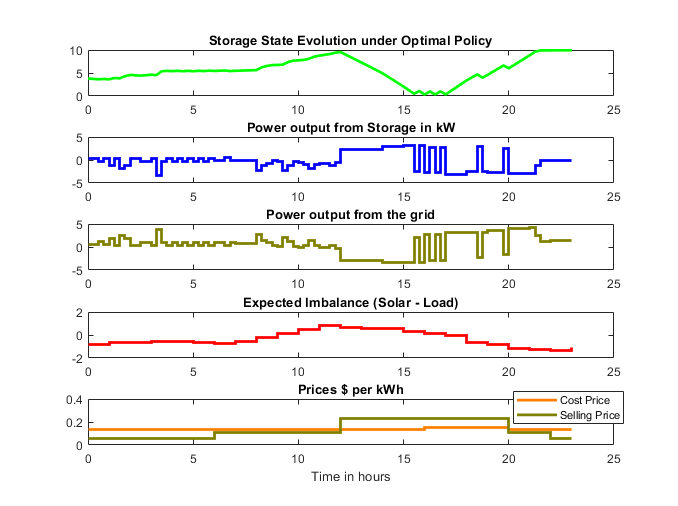}
	\caption{System under the action of near-optimal Policy}
	\label{FIG:optimal}
\end{minipage}\hfill
\begin{minipage}[t]{0.5\textwidth}
	\centering\captionsetup{width=0.9\linewidth}%
	\includegraphics[width=\linewidth]{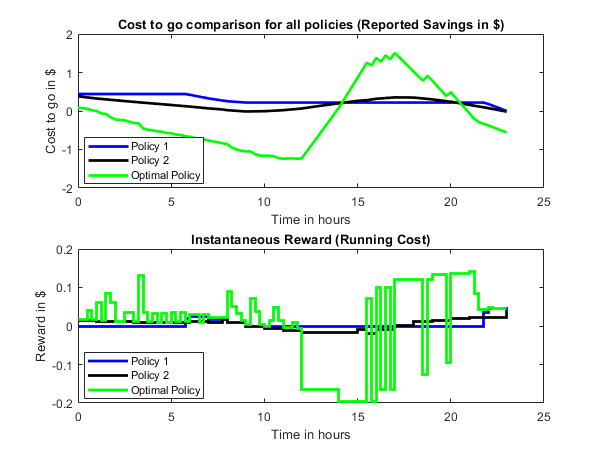}
	\caption{Cost-to-go under various Policies}
	\label{FIG:cost_to_go}
\end{minipage}
\end{figure}

The near-optimal policy described in Section \ref{SEC:OPC} was computed by backward induction. By using the computed values from the next time instant $t_{k+1}$, near-optimal decisions and current values at $t_k$ were computed using the Bellman equation (Equation \ref{EQ:Bellman_eqn}). The near-optimal policy and its action on the state evolution is shown in Figure \ref{FIG:optimal}. It can be observed that the near-optimal policy sells power to the grid when the selling prices are highest, i.e. during the hour interval $[12,20)$. Further, the policy ensures full charge both by the end of the optimization horizon and by noon. The former can be attributed to the large value of terminal cost $g(N,s_N)\ (\mbox{with}\ m=100)$ that would be incurred if the battery were not to be fully charged. However, in case of the latter, the policy charges the battery gradually up to noon despite the usual cost price. This behavior can be understood by noting that the policy gradually charges the battery up to noon in anticipation of maximum profit (minimum cost). This is obtained by selling power to the grid when the selling prices are the highest, which occurs during the hour interval $[12,16)$. Consequently, the near-optimal policy results in a near-minimum cost.  

The cost-to-go associated with the near-optimal policy and the naive policies over the optimization horizon are shown in Figure \ref{FIG:cost_to_go}. It can be observed that the value function in the top subplot (in green) is not necessarily monotonic. This is explained by the reversal of the sign of the grid transaction costs (shown in the bottom subplot). It can also be observed that the near-optimal cost-to-go is greater than its naive counterparts over certain intervals, which might seem contrary to the principle of optimality. However, this is not the case for the following reasons.

First, we note that the value function is a function of the state. Therefore, the values and costs-to-go \emph{at a given instant $t_k \in \mathcal{T}$} must involve the same starting state $s_k$ if a comparison among them is indicative of a comparison among the underlying policies. Clearly, as evident from Figures \ref{fig:policy1}, \ref{fig:policy2}, and \ref{FIG:optimal}, the states at several instants $t_k$ are not necessarily the same except for the instant $t_0$. At $t_0$, a correct comparison of the values $V(0,s_0)$ to the costs-to-go incurred by the naive policies can be made since the same initial battery state $s_0 = 3.8$ were applied to all simulations. In this case, as expected from the principle of optimality the computed near-optimal cost-to-go $V(0,s_0)$ is less than the costs-to-go incurred by the naive policies. Therefore, the costs-to-go under the action of the various policies are consistent with the definition of the value function and the principle of optimality.
\subsubsection{Residential and commercial building simulations over a monthly horizon}
\begin{figure}
\centering
\begin{minipage}[t]{0.5\textwidth}
	\centering\captionsetup{width=0.9\linewidth}%
	\includegraphics[width=\linewidth]{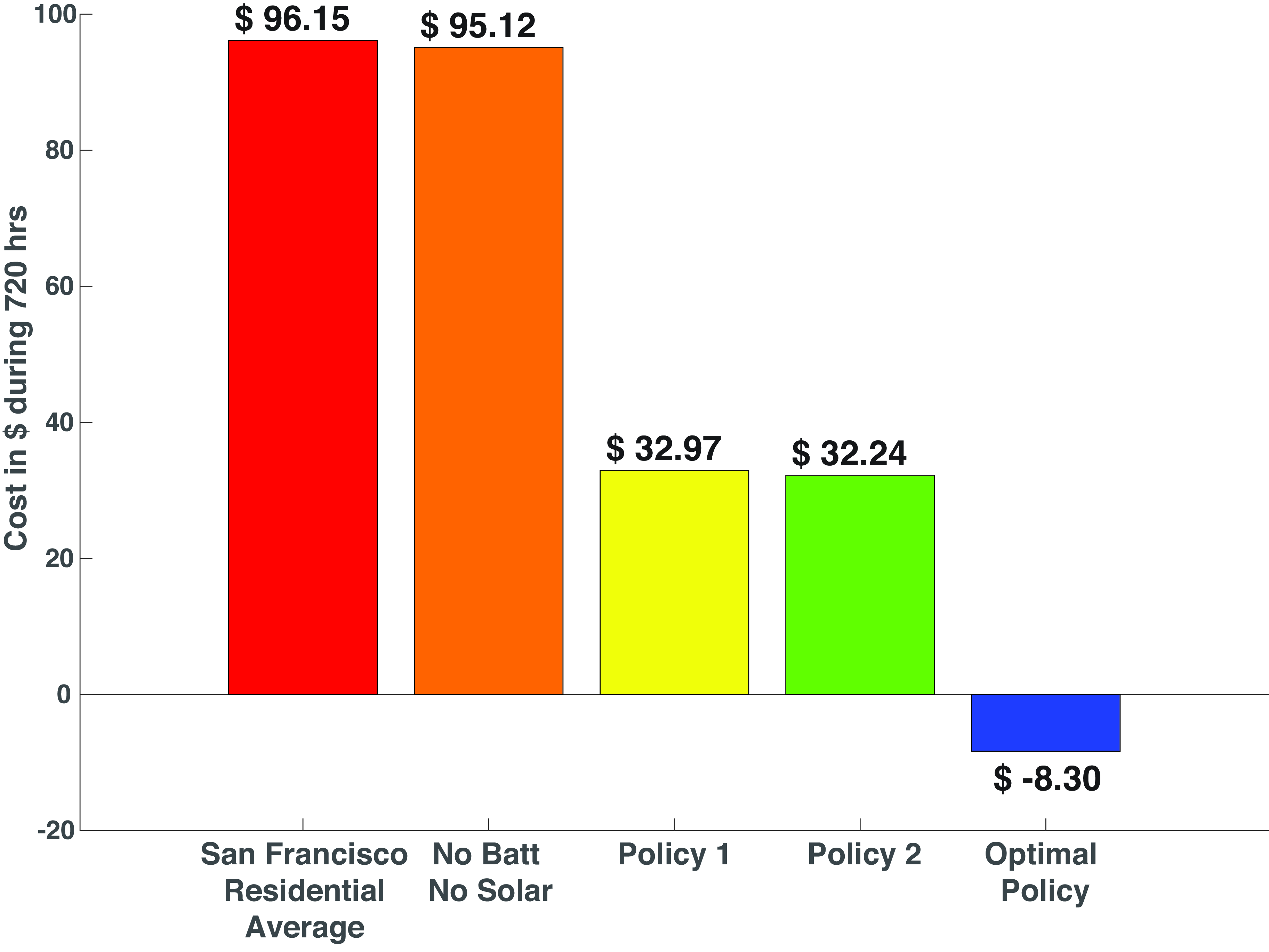}
	\caption{Monthly cost comparison for a residential building}
	\label{fig:barchart}
\end{minipage}\hfill
\begin{minipage}[t]{0.5\textwidth}
	\centering\captionsetup{width=0.9\linewidth}%
	\includegraphics[width=\linewidth]{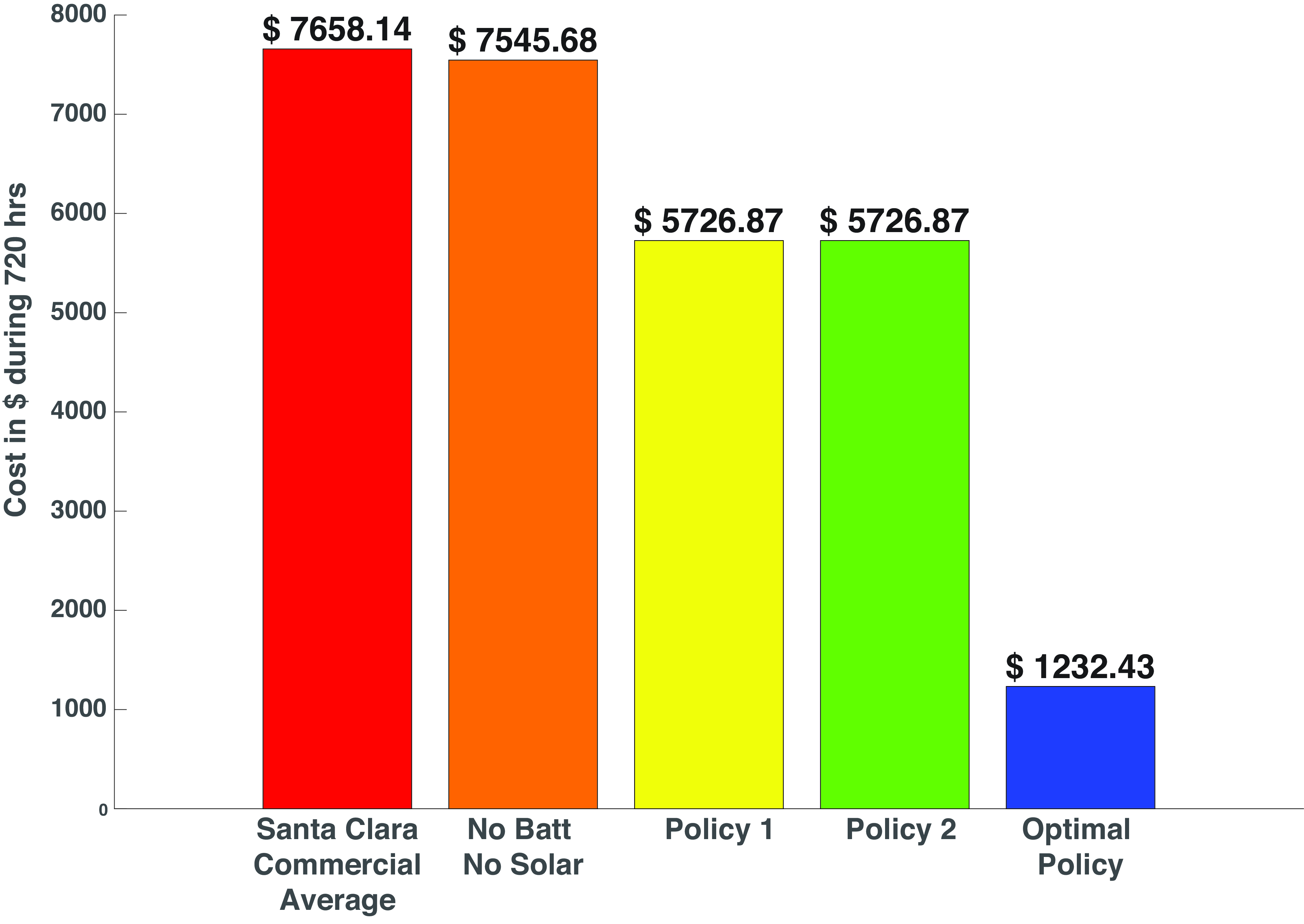}
	\caption{Monthly cost comparison for a commercial building}
	\label{fig:barchart_commercial}
\end{minipage}
\end{figure}

The monthly costs incurred by the various policies in residential and commercial buildings are shown in Figures \ref{fig:barchart} and \ref{fig:barchart_commercial}, respectively. We first compared the \emph{monthly electricity bill from real-world data} (\emph{red bar}) to the \emph{simulated monthly cost} (\emph{orange bar}). For the residential building scenario, the error between San Francisco's average monthly electricity bill ($\$96.15$) \cite{sfo_elec_price} and the corresponding simulated monthly cost ($\$95.12$) was found to be $\approx 1\%$. Similarly, for the commercial building scenario, the error between Santa Clara's average monthly electricity bill ($\$7658.14$) \cite{sc_elec_price} and the corresponding simulated monthly cost ($\$7545.68$) was found to be $\approx 1.5\%$. In these comparisons, the simulation results were based on a nanogrid that is not equipped with either a battery or decision-making ability. In both cases, the comparisons between the experiment baseline (\emph{red bar}) and the simulated data (\emph{orange bar}) indicate a reasonable agreement between the model and real-world data.

In the case of the residential building, the monthly costs under the action of policy 1 (\emph{yellow bar}) and policy 2 (\emph{green bar}) were found to be $\$32.97$ and $\$32.24$, respectively. In addition, we find that the costs-to-go associated with policy 1 (\emph{blue line}) and policy 2 (\emph{black line}) in Figure \ref{FIG:cost_to_go} follow different trajectories but are within close proximity throughout the simulation horizon. This proximity shows that the despite the different considerations of the underlying policies and the different trajectories taken by them, thereby providing further evidence that the overall costs incurred by these policies are expected to be similar under the given demand-supply characteristics. The monthly cost incurred ($\$-8.30$) by the near-optimal policy is shown by the \emph{blue bar} in Figure \ref{fig:barchart}.

Similarly, in case of the commercial building results shown in Figure \ref{fig:barchart_commercial}, the monthly costs incurred by both of the naive policies, shown by the \emph{yellow} and the \emph{green} bars were found to be the same ($\$5726.87$). Further, the costs-to-go for both policies over a daily horizon are shown in top subplot of the Figure \ref{FIG:commercial_heuristic_ctg}, which indicates similarity between the costs incurred and the decisions made by both policies. This can be explained by the other subplots showing a deficit PV generation in the present and expected in the future, thereby resulting in similar decision considerations for both policies. In case of the near-optimal policy shown by the \emph{blue bar}, the cost was found to be $\$1232.43$ which is $\approx 78\%$ lower than the cost incurred by the naive policies.

\subsubsection{Practical utility}
The framework proposed in this work has several types of real-world applications. First, it can be used to solve portfolio optimization problems such as \emph{storage capacity design} and/or \emph{PV capacity design} based on weather uncertainties and pricing schemes. Such solutions can provide recommendations on optimal investment of resources based on Return-On-Investment (ROI) considerations. Second, the MDP framework can be used to recommend an appropriate \emph{utility pricing scheme} that a building must choose in order to minimize expected cost. Third, the cost incurred by the near-optimal policy allows utilities to estimate potential building savings. Such estimates may be used to design \emph{incentive mechanisms} for optimal demand response programming.

%% file: conclusion.tex
\section{Conclusion}
\label{SEC:Conclusion}
The problem of designing optimal policies for nanogrid energy management was treated in the context of stochastic photovoltaic generation, stochastic demand, dynamical storage, and utility pricing. The uncertain photovoltaic generation and the uncertain load were represented by cyclostationary stochastic processes in discrete time and space. The probability distributions of the generation process were inferred from the forecast error data provided by the North American Mesoscale Forecast System. The probability distributions of the load processes were inferred from the yearly demand data provided by OpenEI. Along with a first order dynamical storage model, the demand-supply-storage framework was designed within which decision problems were examined. Naive energy management policies were proposed based on heuristic considerations of demand-supply characteristics. Thereafter, the optimal decision problem was posed and a cost-minimizing policy was sought within a Markov Decision Process framework using stochastic dynamic programming. The stochastic dynamic programming approach was implemented numerically by quantizing the state space, identifying the \emph{determinable feasible decision space}, and approximating the value function, thereby resulting in a near-optimal solution. The near-optimal policy was computed using backward induction. Simulations were carried out for both residential and commercial buildings with real-world parameters for demand, supply, storage, and pricing. Results indicate expected operating cost reductions of $\approx \$40$ and $\approx \$4500$ over a month in residential and commercial building environments, respectively. Future work can investigate the use of continuous state-time formulations, nonlinear battery models, stochastic pricing schemes, partially controllable demand-supply processes, partially observable storage states, grid constraints, and problems pertaining to a network of interconnected nanogrids.

%% file: appendix.tex
\section{Handling constraints in the decision problem}
\label{app:constraint_handling}
The constraints in the decision problem include the following: (i) the power limits $\{P_{min},P_{max}\}$ during the battery charge-discharge process, (ii) the battery state limits $\{0,\mathcal{S}\}$, and (iii) the nanogrid power balance (Equation \ref{EQ:power_balance}). The decisions that respect the above constraints are called \emph{feasible decisions}. Accordingly, the policies resulting in feasible decisions are called \emph{feasible policies}. In order to solve the Bellman equation (Equation \ref{EQ:Bellman_eqn}), a \emph{feasible decision space} is necessary. We construct the \emph{feasible decision space} as follows. At a given time instant $t_k \in \mathcal{T}$, let $e^{(r_e)}_k \in \{e_k^{min},\cdots,e_k^{max}\}$ and $l^{(r_l)}_k \in \{l_k^{min},\cdots,l_k^{max}\}$ represent realizations of the stochastic processes $e_k$ and $l_k$, respectively.
\subsection{Handling State Constraints}
If the state constraint $s_k \in [0,\mathcal{S}]$ is satisfied at $t_k$, we note that the state constraint will be satisfied at $t_{k+1}$ provided the following conditions are satisfied.
\begin{eqnarray}
\label{EQ:SC2CC}
s_{k+1} \in [0,\mathcal{S}]\ \mbox{given}\ s_k \in [0,\mathcal{S}], \notag \\
\iff \eta_s s_k - \xi_p v_k \Delta t \in [0,\mathcal{S}]\ s_k \in [0,\mathcal{S}] \notag \\
(i.e.) -\frac{\eta_s s_k}{\xi_p \Delta t} \leq \ l_k + u_k + e_k \leq \frac{\mathcal{S}-\eta_s s_k}{\xi_p \Delta t} \notag \\
(or) -\frac{\eta_s s_k}{\xi_p \Delta t} - l^{(r_l)}_k - e^{(r_e)}_k \leq \ u_k \leq \frac{\mathcal{S}-\eta_s s_k}{\xi_p \Delta t} - l^{(r_l)}_k - e^{(r_e)}_k
\end{eqnarray}
We write the decision constraint in Equation \ref{EQ:SC2CC} as $u_k \in \Theta(s_k,l^{(r_l)}_k,e^{(r_e)}_k)$. It is easy to observe that, if $u_k \in \Theta(s_k,l^{(r_l)}_k,e^{(r_e)}_k)$ and $s_k \in [0,\mathcal{S}]$, then the above derivation implies that the state constraint is satisfied at $t_{k+1}$.
\subsection{Handling Power Constraints}
By observing that the power flow $v_k$ must be within the limits $[P_{min},P_{max}]$ kW, we can use Equation \ref{EQ:power_balance} to make the following claim:
\begin{eqnarray}
\label{EQ:PC2CC}
v_k \in [P_{min},P_{max}] \notag \\
\iff l^{(r_l)}_k + u_k + e^{(r_e)}_k \in [P_{min},P_{max}] \notag \\
(i.e.)\ P_{min} - l^{(r_l)}_k - e^{(r_e)}_k \leq u_k \leq P_{max} - l^{(r_l)}_k - e^{(r_e)}_k
\end{eqnarray}
Let the decision constraint in Equation \ref{EQ:PC2CC} be written as $u_k \in \Gamma(s_k,l^{(r_l)}_k,e^{(r_e)}_k)$. It is easy to observe that, if $u_k \in \Gamma(s_k,l^{(r_l)}_k,e^{(r_e)}_k)$, then the battery power constraints are satisfied as shown above.

Since the Equations \ref{EQ:SC2CC} and \ref{EQ:PC2CC} constrain the same decision variable $u_k$, the feasible decision space can be obtained as the intersection of these constrained spaces. Let this intersection be represented by $\mathcal{U}_{s_k}^{(r)} := \Theta(s_k,l^{(r_l)}_k,e^{(r_e)}_k) \cap \Gamma(s_k,l^{(r_l)}_k,e^{(r_e)}_k)$. Then $\mathcal{U}_{s_k}^{(r)}$ can be written as,
\begin{equation}
\label{EQ:feasible_decision_space}
max(P_{min},-\frac{\eta_s s_k}{\xi_p \Delta t}) - l^{(r_l)}_k - e^{(r_e)}_k \leq u_k \leq min(P_{max},\frac{\mathcal{S}-\eta_s s_k}{\xi_p \Delta t}) - l^{(r_l)}_k - e^{(r_e)}_k 
\end{equation}
From Equation \ref{EQ:feasible_decision_space}, we observe the following:
\begin{enumerate}
\item The decision space $\ \mathcal{U}_{s_k}^{(r)}$ is guaranteed to have a positive Lebesgue measure, since $max(P_{min},-\frac{\eta_s s_k}{\xi_p \Delta t}) \leq 0$ and $min(P_{max},\frac{\mathcal{S}-\eta_s s_k}{\xi_p \Delta t}) \geq 0$ but both cannot the value 0 simultaneously. Thus, the existence of a non-empty \emph{feasible decision space} is guaranteed by definition.
\item The bounds of the space $\ \mathcal{U}_{s_k}^{(r)}$ depend on the realizations $e^{(r_e)}_k$ and $l^{(r_l)}_k$ of the stochastic processes $e_k$ and $l_k$, respectively. However, during the near-optimal policy design phase, the realizations of the stochastic processes $e^{(r_e)}_k$ and $l^{(r_l)}_k$ are unknown until the time instant $t_k$ occurs in the real world.
\item Despite the guaranteed existence of a \emph{feasible decision space}, it cannot be determined on account of the uncertainty in the PV generation and the load.
\end{enumerate}
In order to eliminate the dependence of the \emph{feasible decision space} on the unknown PV generation and load in the real world at $t_k$, we use the knowledge about the bounds on the corresponding realizations $e^{(r_e)}_k$ and $l^{(r_l)}_k$. Let the range space of the PV generation and the load variables at the time instant $t_k$ be represented by $[e_k^{min},e_k^{max}]$ and $[l_k^{min},l_k^{max}]$, respectively. Since the solar generation and load are bounded in the real world, the bounds $\{e_k^{min},e_k^{max},l_k^{min},l_k^{max}\}$ are physically well-defined. Using these bounds, we construct a subset of the \emph{feasible decision space} $\mathcal{U}_{s_k}^{(r)}$ and denote it as the \emph{determinable feasible decision space} ($\mathcal{U}_{s_k}$).
\begin{equation}
\label{EQ:ControlConstraint}
max(P_{min},-\frac{\eta_s s_k}{\xi_p \Delta t}) - l_k^{min} - e_k^{min} \leq u_k \leq min(P_{max},\frac{\mathcal{S}-\eta_s s_k}{\xi_p \Delta t}) - l_k^{max} - e_k^{max} 
\end{equation}
It is easy to verify that $\mathcal{U}_{s_k}$ is constructed by the intersection of the \emph{feasible decision spaces} across all sample paths with non-zero probability. In other words, $\mathcal{U}_{s_k} := \underset{r=\{0,\cdots,N_r\}}{\cap} \mathcal{U}_{s_k}^{(r)}$, where $\big[P^l(k,l_k=l^{(r_l)})\times P^e(k,e_k=e^{(r_l)})\big] \neq \phi$.

Though the \emph{determinable feasible decision space} $\mathcal{U}_{s_k}$ is a subset of the \emph{feasible decision space} $\mathcal{U}_{s_k}^{(r)}$, it is nevertheless a sufficient condition to ensure that the system constraints are upheld under all possible realizations of the stochastic processes. However, the existence of this space is contingent upon the measure of $\mathcal{U}_{s_k}$ being well-defined. Therefore, the necessary conditions for the existence of a \emph{determinable feasible decision space} with a non-zero measure are:
\begin{eqnarray}
\label{EQ:measure_admissible}
(l_k^{max}+e_k^{max})-(e_k^{min}+l_k^{min}) \notag \\
\leq min(P_{max},\frac{\mathcal{S}-\eta_s s_k}{\xi_p \Delta t}) - max(P_{min},-\frac{\eta_s s_k}{\xi_p \Delta t}) \notag \\
= min(P_{max},\frac{\mathcal{S}-\eta_s s_k}{\xi_p \Delta t}) + min(-P_{min},\frac{\eta_s s_k}{\xi_p \Delta t}) \notag \\
= min(P_{max},\frac{\mathcal{S}-\eta_s s_k}{\xi_p \Delta t}) + min(P_{max},\frac{\eta_s s_k}{\xi_p \Delta t}) \notag \\
\big[\mbox{Assuming}\ P_{min}=-P_{max}\big] \notag \\
= \begin{cases}
    2\times P_{max},& \text{if }  P_{max} < min(\frac{\mathcal{S}-\eta_s s_k}{\xi_p \Delta t},\frac{\eta_s s_k}{\xi_p \Delta t}) \\ \\
   \frac{\mathcal{S}-\eta_s s_k}{\xi_p \Delta t} + P_{max},              & \text{if } \frac{\mathcal{S}-\eta_s s_k}{\xi_p \Delta t} \leq P_{max} < \frac{\eta_s s_k}{\xi_p \Delta t} \\ \\
      P_{max} + \frac{\eta_s s_k}{\xi_p \Delta t},              & \text{if } \frac{\eta_s s_k}{\xi_p \Delta t} \leq P_{max} < \frac{\mathcal{S}-\eta_s s_k}{\xi_p \Delta t} \\ \\
    \frac{\mathcal{S}}{\xi_p \Delta t},              & \text{if } P_{max} > max(\frac{\mathcal{S}-\eta_s s_k}{\xi_p \Delta t},\frac{\eta_s s_k}{\xi_p \Delta t})
\end{cases} \notag \\
= min\Big(2 P_{max},\frac{\mathcal{S}-\eta_s s_k}{\xi_p \Delta t} + P_{max},      P_{max} + \frac{\eta_s s_k}{\xi_p \Delta t},\frac{\mathcal{S}}{\xi_p \Delta t}\Big) \notag \\
\geq min\Big(P_{max},\frac{\mathcal{S}}{\xi_p \Delta t}\Big) \\
\big[\mbox{equality iff}\ s_k = 0\ or\ \eta_s s_k = \mathcal{S}\ or\ \frac{\mathcal{S}}{\xi_p \Delta t} \leq P_{max}\big] \notag
\end{eqnarray}
It is easy to verify that, if $(l_k^{max}+e_k^{max})-(e_k^{min}+l_k^{min}) \leq min(P_{max},\frac{\mathcal{S}}{\xi_p \Delta t})$, then the inequation \ref{EQ:ControlConstraint} holds true. Hence, the inequation $(l_k^{max}+e_k^{max})-(e_k^{min}+l_k^{min}) \leq min(P_{max},\frac{\mathcal{S}}{\xi_p \Delta t})$ provides a stronger condition for the existence of a non-empty \emph{determinable feasible decision space} $\mathcal{U}_{s_k}$. Progressively stronger sufficiency conditions can be derived as follows:
\begin{eqnarray}
\label{EQ:suff_cond}
(l_k^{max}+e_k^{max})-(e_k^{min}+l_k^{min}) \leq min\Big(P_{max},\frac{\mathcal{S}}{\xi_p \Delta t}\Big),\ or\ \\
\label{EQ:suff_cond2}
e_k^{max}-l_k^{min} \leq min\Big(P_{max},\frac{\mathcal{S}}{\xi_p \Delta t}\Big) \\
\big[\mbox{since}\ -e_k^{min}\leq0, l_k^{max}\leq0\ \mbox{by definition}\big] \notag \\
\label{EQ:suff_cond3}
(i.e.)\ \max_{k \in \{0,\cdots,N\}}(e_k^{max}-l_k^{min}) \leq min\Big(P_{max},\frac{\mathcal{S}}{\xi_p \Delta t}\Big) \notag \\
\big[\mbox{Let}\ \bar{l}_k=-l_k,\ \mbox{and let}\ \bar{l}_k \in \{\bar{l}^{min},\cdots,\bar{l}^{max}\}\big] \notag \\
\max_{k \in \{0,\cdots,N\}} (e_k^{max}+\bar{l}_k^{max}) \leq min\Big(P_{max},\frac{\mathcal{S}}{\xi_p \Delta t}\Big) \notag \\
\big[\mbox{Since}\ \max_{k \in \{0,\cdots,N\}} \bar{l}_k^{max} = -l^{max}\big] \notag \\
(e^{max}-l^{min}) \leq min\Big(P_{max},\frac{\mathcal{S}}{\xi_p \Delta t}\Big) 
\end{eqnarray}
Note that the left-hand side (LHS) of the Equation \ref{EQ:suff_cond} represents the maximum possible demand-supply offset gap at the instant $t_k \in \mathcal{T}$. On the other hand, the LHS of the Equation \ref{EQ:suff_cond2} represents the maximum demand-supply sum at the instant $t_k$. The LHS of the strongest condition derived in Equation \ref{EQ:suff_cond3} represents the maximum possible demand-supply sum over the entire partition $\mathcal{T}$. In each of these sufficiency conditions, it is imperative to note that the right-hand side represents a time-independent subset of the storage and the duration parameters. Let this subset consisting of $\{P_{max},\mathcal{S},\xi_p,\Delta t\}$ be denoted by $\lambda$.

Equations \ref{EQ:ControlConstraint}-\ref{EQ:suff_cond3} represent the worst case sufficiency conditions for the existence of a corresponding \emph{determinable feasible decision space} $\mathcal{U}_{s_k}$. Further, for any grid transaction decision $u_k$ within $\mathcal{U}_{s_k}$, it is guaranteed that the battery power constraint is satisfied, since $\mathcal{V}_{s_k} = [P_{min},P_{max}]$ is equivalent to Equation \ref{EQ:PC2CC} by definition. 

Let the parameters $\lambda$ satisfying the sufficiency condition\footnote{Since the conditions represented by the Equations \ref{EQ:ControlConstraint}-\ref{EQ:suff_cond3} are progressively stronger, satisfying any of these conditions ensures that the Equation \ref{EQ:ControlConstraint} is satisfied.} in Equation \ref{EQ:suff_cond3} belong to the space $\Lambda \subseteq \mathbb{R}^{dim(\lambda)}$, where $dim(\lambda)$ refers to the dimension of $\lambda$. We refer to $\Lambda$ as the \emph{determinable feasible configuration space} of the system. Therefore, the presence of the system parameters within the \emph{determinable feasible configuration space} ($\lambda \in \Lambda$) guarantees the existence of a corresponding \emph{determinable feasible decision space} which is required to design the near-optimal policy. Only decisions belonging to the \emph{determinable feasible decision space} are considered admissible for computing the near-optimal policy, hence we also refer to this decision space $(\mathcal{U}_{s_k},\mathcal{V}_{s_k})$ as the \emph{admissible decision space} for the optimal energy management problem.

In summary, a \emph{feasible decision space} is guaranteed to exist for every realization of the stochastic processes $e_k$ and $l_k$ (see Equation \ref{EQ:feasible_decision_space}). However, such realizations are unknown during the design of the near-optimal policy. This is because such design involves computing the expected state at the next instant of time $\mathbb{E}(\hat{s}_{k+1})$ in the Bellman Equation (Equation \ref{EQ:Bellman_eqn}). Thus, every possible realization of the load and PV stochastic processes is considered. Accordingly, every resulting realization of $\hat{s}_{k+1}$ (= $s_{k+1}$) is expressed as a function of the decision $u_k$ given the knowledge of $s_k$. By ensuring that every possible realization of $\hat{s}_{k+1}$ is expressed as a function of $u_k$, and that $u_k$ satisfies the constraints in Equation \ref{EQ:SC2CC} and Equation \ref{EQ:PC2CC}, the state constraints are translated into corresponding control constraints $u_k$. In this manner, sufficiency conditions for satisfying the system constraints are derived. Finally, we note that computing a near-optimal policy using SDP is feasible only if the conditions in Equations \ref{EQ:measure_admissible}-\ref{EQ:suff_cond3} are satisfied, thereby ensuring a non-empty \emph{admissible decision space} ($\mathcal{U}_{s_k} \neq \phi$).